# Boundary-to-Displacement Asymptotic Gains for Wave Systems With Kelvin-Voigt Damping

**Iasson Karafyllis[*], Maria Kontorinaki[**] and Miroslav Krstic[***]**

[*]Dept. of Mathematics, National Technical University of Athens, Zografou Campus, 15780, Athens, Greece, email: iasonkar@central.ntua.gr

[**]Department of Statistics and Operations Research, University of Malta, Tal-Qroqq Campus, Msida, 2080, Malta, email: kontorinmaria@gmail.com

[***]Dept. of Mechanical and Aerospace Eng., University of California, San Diego, La Jolla, CA 92093-0411, U.S.A., email: krstic@ucsd.edu

## Abstract

We provide estimates for the asymptotic gains of the displacement of a vibrating string with endpoint forcing, modeled by the wave equation with Kelvin-Voigt and viscous damping and a boundary disturbance. Two asymptotic gains are studied: the gain in the $L^2$ spatial norm and the gain in the spatial sup norm. It is shown that the asymptotic gain property holds in the $L^2$ norm of the displacement without any assumption for the damping coefficients. The derivation of the upper bounds for the asymptotic gains is performed by either employing an eigenfunction expansion methodology or by means of a small-gain argument, whereas a novel frequency analysis methodology is employed for the derivation of the lower bounds for the asymptotic gains. The graphical illustration of the upper and lower bounds for the gains shows that that the asymptotic gain in the $L^2$ norm is estimated much more accurately than the asymptotic gain in the sup norm.

## 1. Introduction

Asymptotic gain properties for finite-dimensional systems were introduced in [1,3,14,21]. More specifically, in [21] it was shown the Input-to-State Stability (ISS) superposition theorem, which was extended in [1] for the case of the output asymptotic gain property. Recently, the asymptotic gain property has been used in time-delay systems (see [9]) and abstract infinite-dimensional systems (see [14]). For linear systems, where the asymptotic gain function is linear, the asymptotic gain (coefficient) may be used as a measure of the sensitivity of the system with respect to external disturbances.

The wave equation with viscous and Kelvin-Voigt damping is the prototype Partial Differential Equation (PDE) for the description of vibrations in elastic media with energy dissipation but may also arise in different physical problems (e.g., movement of chemicals underground; see [5,10]). The study of the dynamics of the wave equation with viscous and Kelvin-Voigt damping has attracted the interest of many researchers (see [2,4,6,13,15,16,22]). Recently, the control of the wave equation with Kelvin-Voigt damping was studied in [12,20], while the control of the wave equation with viscous damping was studied in [17,18,19].

The long-time behavior of the wave equation with viscous and Kelvin-Voigt damping under a boundary disturbance was studied in [10], within the theoretical framework of the ISS property (see [7,11]). Moreover, an upper bound of the maximum displacement that can be caused by a unit boundary disturbance (Input-to-Output Stability gain in the sup norm) was given under a specific assumption for the damping coefficients.

The aim of the present work is the extension of the result in [10] to various directions by employing the asymptotic gain property for the wave equation with viscous and Kelvin-Voigt damping given by

$$\frac{\partial^2 u}{\partial t^2}(t,x) = \frac{\partial^2 u}{\partial x^2}(t,x) + \sigma \frac{\partial^3 u}{\partial t \partial x^2}(t,x) - \mu \frac{\partial u}{\partial t}(t,x), \text{ for } (t,x) \in (0,+\infty) \times (0,1),  \quad (1)$$

$$\begin{matrix} u(t,0) = d(t) \\ u(t,1) = 0 \end{matrix}, \text{ for } t \geq 0, \quad (2)$$

where $\sigma > 0$, $\mu \geq 0$ are constants. This is the mathematical model of a vibrating string with internal and external damping and $u(t,x)$ denotes the displacement of the string at time $t \geq 0$ and position $x \in [0,1]$. One end of the string (at $x=1$) is pinned down while external forcing acts on the other end of the string (at $x=0$). The effect of the external forcing is described by the boundary disturbance $d(t)$ and

- we show that the asymptotic gain property holds in the $L^2$ norm of the displacement without any assumption for the damping coefficients (Theorem 3.1),
- we provide upper and lower bounds for the asymptotic gains in the sup-norm and $L^2$ norm of the displacement (Theorems 2.3 and 2.4 and Theorems 3.1 and 3.2).

*This is the first time that a systematic study for the asymptotic gains of a PDE with a boundary disturbance is being performed.* For example, Fig. 1 and Fig. 2 depict the lower and upper bounds of the asymptotic gain, in the sup-norm ($L_\infty$ and $U_\infty$) and $L^2$ norm ($L_2$ and $U_2$), respectively, for $\sigma = 1$ and for a wide range of values for $\mu$.

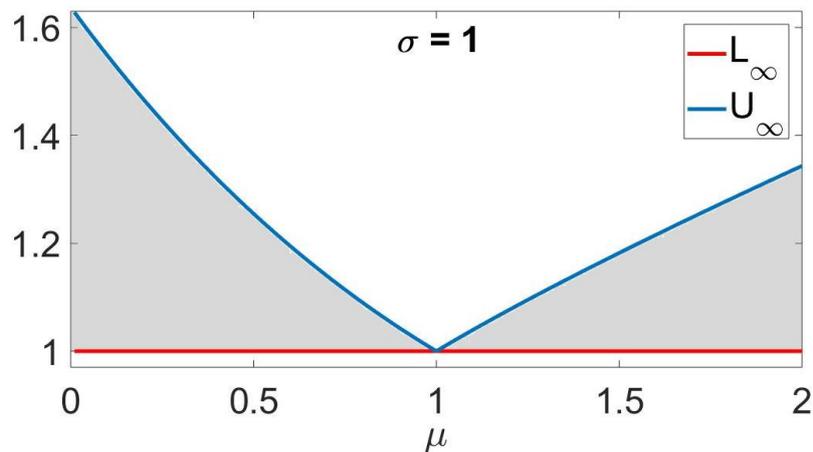

**Fig. 1:** The lower (red) and upper (blue) bounds for the asymptotic gain in the spatial sup-norm for $\sigma = 1$. The grey area depicts possible values for the asymptotic gain in the sup-norm.

The derivation of the upper bounds for the asymptotic gains is performed by either employing an eigenfunction expansion methodology (see also [8,11]) or by means of a small-gain argument (as in [10,11]). The derivation of the lower bounds for the asymptotic gains is performed by means of a novel methodology, *which has never used before* for boundary disturbances in PDEs: a frequency analysis methodology. The methodology is convenient and can be used for other PDEs.

The structure of the present work is as follows: Section 2 introduces the reader to the notion of the asymptotic gain and describes the frequency analysis methodology. All main results in the sup-norm are stated in Section 2, while all main results in the $L^2$ norm are given in Section 3. The proofs of all main results are provided in Section 4. Section 5 provides a graphical illustration of the main results provided in Sections 2 and 3. Finally, the concluding remarks of the present work are contained in Section 6.

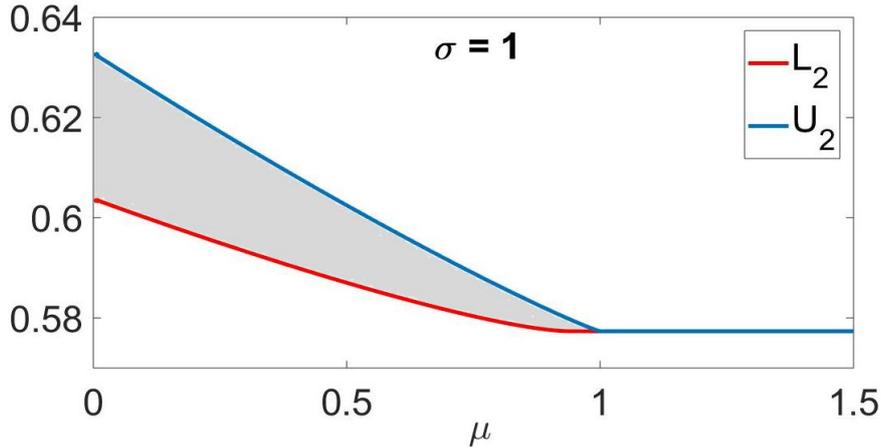

**Fig. 2:** The lower (red) and upper (blue) bounds for the asymptotic gain in the spatial $L^2$ norm for $\sigma = 1$. The grey area depicts possible values for the asymptotic gain in the $L^2$ norm.

**Notation.** Throughout this paper, we adopt the following notation.

* $\Re_+ := [0, +\infty)$. For every $x \in \Re$, $[x]$ denotes the integer part of $x$.
* Let $A \subseteq \Re^n$ be an open set and let $U \subseteq \Re^n$, $\Omega \subseteq \Re$ be sets with $A \subseteq U \subseteq \bar{A}$. By $C^0(U)$ (or $C^0(U;\Omega)$), we denote the class of continuous mappings on $U$ (which take values in $\Omega$). By $C^k(U)$ (or $C^k(U;\Omega)$), where $k \geq 1$, we denote the class of continuous functions on $U$, which have continuous derivatives of order $k$ on $U$ (and also take values in $\Omega$).
* $L^2(0,1)$ denotes the equivalence class of measurable functions $f:[0,1] \to \Re$ for which $\|f\|_2 = \left(\int_0^1 |f(x)|^2 dx\right)^{1/2} < +\infty$. $L^\infty(0,1)$ denotes the equivalence class of measurable functions $f:[0,1] \to \Re$ for which $\|f\|_\infty := \operatorname{ess\,sup}_{x \in (0,1)}(|f(x)|) < +\infty$.
* Let $u: \Re_+ \times [0,1] \to \Re$ be given. We use the notation $u[t]$ to denote the profile at certain $t \geq 0$, i.e. $(u[t])(x) = u(t,x)$, for all $x \in [0,1]$. When $u(t,x)$ is differentiable with respect to $x \in [0,1]$, we use the notation $u'(t,x)$ for the derivative of $u$ with respect to $x \in [0,1]$, i.e., $u'(t,x) = (\partial u/\partial x)(t,x)$.
* For an integer $k \geq 1$, $H^k(0,1)$ denotes the Sobolev space of functions in $L^2(0,1)$ with all its weak derivatives up to order $k \geq 1$ in $L^2(0,1)$.

## 2. Asymptotic Gain in the Sup Norm

*2.1. Background*

Let (P) be a specific 1-D PDE control system with domain $[0,1]$ and let $d(t) \in \Re^m$ be a disturbance input that enters the system. Let $X$ be the state space of the system ($X$ is a normed linear functional space of functions defined on $(0,1)$) and for each $(u_0, w_0) \in X$ let $\Phi(u_0, w_0) \neq \emptyset$ denote the set of bounded disturbances $d: \Re_+ \to \Re^m$ for which the system has a unique solution $(u[t], w[t]) \in X$ for $t \geq 0$ (a continuous mapping $\Re_+ \ni t \to (u[t], w[t]) \in X$ with $(u[0], w[0]) = (u_0, w_0)$ for which the classical semigroup property holds). We assume that there exists at least one initial state $(u_0, w_0) \in X$ and a non-zero bounded disturbance $d \in \Phi(u_0, w_0)$ for which the system has a unique solution $(u[t], w[t]) \in X$ for $t \geq 0$ (because otherwise the only allowable disturbance would be the zero input).

Suppose that (P) satisfies an output asymptotic gain property with linear gain function, i.e., suppose that there exists a constant $\gamma \geq 0$ such that for every initial state $(u_0, w_0) \in X$ and for every allowable disturbance $d \in \Phi(u_0, w_0)$ the following inequality holds for the solution $(u[t], w[t]) \in X$:

$$\limsup_{t \to +\infty} \left( \|u[t]\| \right) \leq \gamma \sup_{t \geq 0} \left( |d(t)| \right), \tag{3}$$

where $\|u[t]\|$ is a norm of a functional space of functions defined on $(0,1)$. Two things should be emphasized at this point:

1) Property (3) is an *output* asymptotic gain property. Notice that the state may have additional components, denoted here by $w[t]$. The output map is the map $X \ni (u, w) \to u \in S$, where $S$ is a normed linear functional space (of functions defined on $(0,1)$) with norm $\|\ \|$. We assume that the map $X \ni (u, w) \to u \in S$ is continuous. Notice that even in the case where the only component of the state is $u[t]$, the output map is not the identity mapping when norm $\|\ \|$ does not coincide with the norm of $X$.

2) Property (3) is an output asymptotic gain property *with linear gain function*. The function $g\left(\sup_{t \geq 0}(|d(t)|)\right) = \gamma \sup_{t \geq 0}(|d(t)|)$ that appears on the right hand side of (3) is linear. This feature allows us to use the name "$d$-to-$u$ asymptotic gain in the norm of $S$" for the quantity

$$\gamma_{\|\ \|} := \sup \left\{ \frac{\limsup_{t \to +\infty}(\|u[t]\|)}{\sup_{t \geq 0}(|d(t)|)} : (u_0, w_0) \in X, d \in \Phi(u_0, w_0), d \neq 0 \right\} \leq \gamma, \tag{4}$$

where $(u[t], w[t]) \in X$ denotes the solution with $(u[0], w[0]) = (u_0, w_0)$, corresponding to $d \in \Phi(u_0, w_0)$.

It should be noticed that $\gamma_{\|\ \|}$, i.e., the $d$-to-$u$ asymptotic gain in the norm of $S$ defined by (4), is the smallest constant $\gamma \geq 0$ for which (3) holds for every initial state $(u_0, w_0) \in X$ and for every bounded disturbance $d: \Re_+ \to \Re^m$ for which the system has a unique solution $(u[t], w[t]) \in X$ for $t \geq 0$.

**Remark 2.1:** Due to the semigroup property, the output asymptotic gain property with linear gain function (3) holds if and only if the following inequality holds for every initial state $(u_0, w_0) \in X$ and for every $d \in \Phi(u_0, w_0)$:

$$\limsup_{t \to +\infty}(\|u[t]\|) \leq \gamma \limsup_{t \to +\infty}(|d(t)|). \tag{5}$$

Therefore, the $d$-to-$u$ asymptotic gain in the norm of $S$ defined by (4) satisfies the following equation

$$\gamma_{\|\ \|} = \sup\left\{ \frac{\limsup_{t \to +\infty}(\|u[t]\|)}{\limsup_{t \to +\infty}(|d(t)|)} : (u_0, w_0) \in X, d \in \Phi(u_0, w_0), \limsup_{t \to +\infty}(|d(t)|) > 0 \right\}$$

provided that there exists at least one initial state $(u_0, w_0) \in X$ and a disturbance $d \in \Phi(u_0, w_0)$ with $\limsup_{t \to +\infty}(|d(t)|) > 0$. If for every $(u_0, w_0) \in X$ the allowable disturbance set $\Phi(u_0, w_0)$ contains only disturbances $d : \Re_+ \to \Re^m$ with $\limsup_{t \to +\infty}(|d(t)|) = 0$, then definition (4) and inequality (5) imply that $\gamma_{\|\ \|} = 0$.

*2.2. Frequency Analysis*

The frequency analysis of the $d$-to-$u$ asymptotic gain in the norm of $S$ consists of two steps.

<u>Step 1:</u> Let $d_\omega : \Re_+ \to \Re^m$ be a parameterized family of non-zero inputs, with parameter $\omega > 0$, which are periodic with period $T = 2\pi/\omega$, i.e., $d_\omega(t + 2\pi/\omega) = d_\omega(t)$ for all $t \geq 0$. For each $\omega > 0$, find a periodic solution $(u_\omega[t], w_\omega[t]) \in X$ of (P) with period $T = 2\pi/\omega$ that corresponds to the input $d_\omega : \Re_+ \to \Re^m$.

<u>Step 2:</u> If Step 1 can be accomplished, then estimate the $d$-to-$u$ asymptotic gain in the norm of $S$, by means of the following inequality

$$\sup_{\omega > 0}\left( \frac{\max_{0 \leq t \leq 2\pi/\omega}(\|u_\omega[t]\|)}{\sup_{0 \leq t < 2\pi/\omega}(|d_\omega(t)|)} \right) \leq \gamma_{\|\ \|}. \tag{6}$$

Inequality (6) is a direct consequence of (4), the fact that $(u_\omega[t], w_\omega[t]) \in X$ is a periodic solution of (P) with period $T = 2\pi/\omega$ that corresponds to the non-zero input $d_\omega : \Re_+ \to \Re^m$ and the fact that the mappings $X \ni (u, w) \to u \in S$, $\Re_+ \ni t \to (u[t], w[t]) \in X$ are continuous mappings (and therefore $\max_{0 \leq t \leq 2\pi/\omega}(\|u_\omega[t]\|)$ exists for each $\omega > 0$).

*2.3. Wave Equation with Kelvin-Voigt and Viscous Damping*

Consider the wave equation with Kelvin-Voigt and viscous damping (1), (2), where $\sigma > 0$, $\mu \geq 0$ are constants. Without loss of generality, the tension, i.e., the coefficient of $(\partial^2 u/\partial x^2)(t, x)$ in (1) has been assumed to be equal to 1; this can always be achieved with appropriate time scaling.

In order to obtain an existence/uniqueness result for the wave equation with Kelvin-Voigt and viscous damping we first need to move the disturbance from the boundary to the domain and make the boundary conditions homogeneous. However, this process should be done with caution because we would also like the non-homogeneous term that will appear in the PDE to be expressed by a sufficiently fast convergent Fourier series. This may be done in the following way: if we assume that $d \in C^4(\Re_+)$ then we may define the functions

$$\beta_{i+1}(t) = \sigma^2 \dot{\beta}_i(t) + \sigma(\mu\sigma - 1)\beta_i(t) - (\mu\sigma - 1)\int_0^t \exp\left(-\frac{t-s}{\sigma}\right)\beta_i(s)ds, \text{ for } i = 0,1, \ t \geq 0, \quad (7)$$

with

$$\beta_0(t) = -\mu\dot{d}(t) - \ddot{d}(t), \text{ for } t \geq 0, \quad (8)$$

and perform the following transformation

$$u(t,x) = w(t,x) + d(t)(1-x) + \frac{1}{6\sigma}(x^3 - 3x^2 + 2x)\int_0^t \exp\left(-\frac{t-s}{\sigma}\right)\beta_0(s)ds$$

$$-\frac{1}{360\sigma^4}(8x - 20x^3 + 15x^4 - 3x^5)\int_0^t \exp\left(-\frac{t-s}{\sigma}\right)\beta_1(s)ds \quad (9)$$

which transforms problem (1), (2) to the following problem

$$\frac{\partial^2 w}{\partial t^2}(t,x) = \frac{\partial^2 w}{\partial x^2}(t,x) + \sigma\frac{\partial^3 w}{\partial t \partial x^2}(t,x) - \mu\frac{\partial w}{\partial t}(t,x) + \frac{8x - 20x^3 + 15x^4 - 3x^5}{360\sigma^6}\beta_2(t),$$

$$\text{for } (t,x) \in (0,+\infty) \times (0,1), \quad (10)$$

$$w(t,0) = w(t,1) = 0, \text{ for } t \geq 0. \quad (11)$$

For every pair of functions $f \in H^6(0,1)$, $g \in H^4(0,1)$ with $f(0) = f''(0) = f^{(4)}(0) = 0$, $f(1) = f''(1) = f^{(4)}(1) = 0$, $g(0) = g''(0) = 0$, $g(1) = g''(1) = 0$, a solution $w \in C^2(\Re_+ \times [0,1])$ with $(\partial w/\partial t)[t] \in C^2([0,1])$ for all $t \geq 0$ of the initial-boundary value problem (10), (11) with

$$w[0] = f, \quad \frac{\partial w}{\partial t}[0] = g \quad (12)$$

can be expressed by the formula

$$w(t,x) = \sqrt{2}\sum_{n=1}^{\infty} a_n(t)\sin(n\pi x), \text{ for } t \geq 0, x \in [0,1] \quad (13)$$

where $a_n(t)$ for $n = 1,2,...$ are the solutions of

$$\ddot{a}_n(t) + (\mu + n^2\pi^2\sigma)\dot{a}_n(t) + n^2\pi^2 a_n(t) = \frac{\beta_2(t)\sqrt{2}}{360\sigma^6}\int_0^1 (8x - 20x^3 + 15x^4 - 3x^5)\sin(n\pi x)dx \quad (14)$$

with

$$a_n(0) = \sqrt{2}\int_0^1 f(x)\sin(n\pi x)dx, \quad \dot{a}_n(0) = \sqrt{2}\int_0^1 g(x)\sin(n\pi x)dx, \text{ for } n = 1,2,.... \quad (15)$$

Indeed, since $a_n(0) = O(n^{-6})$ (a consequence of the fact that $f \in H^6(0,1)$ with $f(0) = f''(0) = f^{(4)}(0) = 0$, $f(1) = f''(1) = f^{(4)}(1) = 0$), $\dot{a}_n(0) = O(n^{-4})$ (a consequence of the fact that $g \in H^4(0,1)$ with $g(0) = g''(0) = 0$, $g(1) = g''(1) = 0$) and since $c_n := (\sqrt{2}/360\sigma^6) \cdot \int_0^1 (8x - 20x^3 + 15x^4 - 3x^5)\sin(n\pi x)dx = O(n^{-5})$, it can be shown that the solution of (14), (15) satisfies $a_n(t) = O(n^{-6})$ for every fixed $t \geq 0$. This follows from the fact that for sufficiently large $n \geq 1$ (so that $\mu + n^2\pi^2\sigma > 2n\pi$), we have for $t \geq 0$:

$$a_n(t) = \frac{a_n(0)}{2r_n}\left((k_n + r_n)\exp(-(k_n - r_n)t) + (r_n - k_n)\exp(-(k_n + r_n)t)\right)$$

$$+\frac{\dot{a}_n(0)}{2r_n}\left(\exp(-(k_n - r_n)t) - \exp(-(k_n + r_n)t)\right) + \frac{c_n}{2r_n}\int_0^t \left(\exp(-(k_n - r_n)(t-s)) - \exp(-(k_n + r_n)(t-s))\right)\beta_2(s)ds$$

where $k_n = (\mu + n^2\pi^2\sigma)/2$, $r_n = \sqrt{(\mu + n^2\pi^2\sigma)^2 - 4n^2\pi^2}/2$ (notice that $c_n/r_n = O(n^{-7})$ and that there exists a constant $L > 0$ so that $k_n - r_n \geq L$ for all sufficiently large $n \geq 1$). Therefore, the function $w$ defined by

(13) is of class $C^2(\Re_+ \times [0,1])$ with $(\partial w/\partial t)[t] \in C^2([0,1])$ for all $t \geq 0$ and satisfies (10), (11) for all $t \geq 0$, $x \in [0,1]$.

Using the transformation (9), we are in a position to give the following proposition.

**Proposition 2.2:** *Suppose that* $d \in C^4(\Re_+)$ *and* $f \in H^6(0,1)$, $g \in H^4(0,1)$ *with* $f(0) = d(0)$, $f''(0) = f^{(4)}(0) = 0$, $f(1) = f''(1) = f^{(4)}(1) = 0$, $g(0) = \dot{d}(0)$, $g''(0) = (1/\sigma)(\mu \dot{d}(0) + \ddot{d}(0))$, $g(1) = g''(1) = 0$. *Then the initial-boundary value problem (1), (2) with*

$$u[0] = f, \quad \frac{\partial u}{\partial t}[0] = g \tag{16}$$

*has a unique solution* $u \in C^2(\Re_+ \times [0,1])$ *with* $(\partial u/\partial t)[t] \in C^2([0,1])$ *for all* $t \geq 0$, *which satisfies (1), (2) for all* $t \geq 0$, $x \in [0,1]$.

**Proof:** A straightforward application of transformation (9) and formulas (7), (8). Uniqueness follows by considering the energy functional

$$E(t) = \int_0^1 \left( \left( \frac{\partial u}{\partial t}(t,x) - \frac{\partial v}{\partial t}(t,x) \right)^2 + \left( \frac{\partial u}{\partial x}(t,x) - \frac{\partial v}{\partial x}(t,x) \right)^2 \right) dx$$

where $u, v \in C^2(\Re_+ \times [0,1])$ with $(\partial u/\partial t)[t] \in C^2([0,1])$, $(\partial v/\partial t)[t] \in C^2([0,1])$ for all $t \geq 0$, are two arbitrary solutions the initial-boundary value problem (1), (2), (16). Since $\dot{E}(t) \leq 0$ for all $t \geq 0$ and since $E(0) = 0$, we obtain that $E(t) = 0$ for all $t \geq 0$. This implies that $(\partial u/\partial x)(t,x) = (\partial v/\partial x)(t,x)$ for all $t \geq 0$, $x \in [0,1]$, which gives $u \equiv v$. The proof is complete. ◁

Let $X$ be the linear space of all $(f,g) \in C^2([0,1]) \times C^2([0,1])$ for which there exists a bounded input $d \in C^2(\Re_+)$ such that the initial-boundary value problem (1), (2), (16) has a unique solution $u \in C^2(\Re_+ \times [0,1])$ with $(\partial u/\partial t)[t] \in C^2([0,1])$ for all $t \geq 0$. Notice that Proposition 2.1 guarantees that $D \subseteq X$, where

$$D = \{(f,g) \in H^6(0,1) \times H^4(0,1) : f''(0) = f^{(4)}(0) = f(1) = f''(1) = f^{(4)}(1) = g(1) = g''(1) = 0\}.$$

For every $(f,g) \in X$, let $\Phi(f,g)$ denote the set of bounded disturbances $d \in C^2(\Re_+)$ for which the initial-boundary value problem (1), (2), (16) has a unique solution $u \in C^2(\Re_+ \times [0,1])$ with $(\partial u/\partial t)[t] \in C^2([0,1])$ for all $t \geq 0$.

*2.4. Asymptotic Gain in the sup norm*

We assume that the following assumption holds.

**(H1)** *There exists a constant* $\gamma \geq 0$ *such that for every bounded disturbance* $d \in C^2(\Re_+)$ *for which (1), (2) has a unique solution* $u \in C^2(\Re_+ \times [0,1])$ *with* $(\partial u/\partial t)[t] \in C^2([0,1])$ *for all* $t \geq 0$, *the following inequality holds:*

$$\limsup_{t \to +\infty} (\|u[t]\|_\infty) \leq \gamma \sup_{t \geq 0} (|d(t)|). \tag{17}$$

It was shown in [10] that Assumption (H1) holds when $2 < 2\mu\sigma + \sigma^2 \pi^2$. When Assumption (H1) holds, we are in a position to define the asymptotic gain of the displacement in the sup norm by means of the formula

$$\gamma_\infty := \sup\left\{ \frac{\limsup_{t\to+\infty}\left(\|u[t]\|_\infty\right)}{\sup_{t\geq 0}\left(|d(t)|\right)} : (f,g)\in X, d\in \Phi(f,g), d\neq 0 \right\}, \tag{18}$$

where $u\in C^2(\Re_+\times[0,1])$ with $(\partial u/\partial t)[t]\in C^2([0,1])$ for all $t\geq 0$ is the solution of the initial-boundary value problem (1), (2), (16). The results in [10] show that when $2<2\mu\sigma+\sigma^2\pi^2$, the following inequality holds

$$\gamma_\infty \leq g\left(\frac{\mu\sigma-1}{\sigma^2}\right), \tag{19}$$

where $g(s):=\inf\left\{1\big/\left(\sin(\theta)\left(1-\sqrt{p(s,\theta)}\right)^2\right):0<\theta<\left(\pi-\sqrt{|s|-s}\right)/2\right\}$ and $p(s,\theta):=|s|\big/\left(s+(\pi-2\theta)^2\right)$.

Our first main result in the present work concerning the sup norm is a sharper result than the result given in [10]. Its proof is provided in Section 4.

**Theorem 2.3:** *Consider the wave equation with Kelvin-Voigt and viscous damping (1), (2), where $\sigma>0$, $\mu\geq 0$ are constants with $2<2\mu\sigma+\sigma^2\pi^2$. Then Assumption (H1) holds and inequality (19) holds with $g(s):=\inf\left\{1\big/\left(\sin(\theta)(1-p(s,\theta))\right):0<\theta<\pi-\sqrt{|s|-s}\right\}$ and $p(s,\theta):=|s|\big/\left(s+(\pi-\theta)^2\right)$.*

Our second main result in the present work concerning the sup norm is given below and it is proved by means of the frequency analysis procedure that was described above. Its proof is provided in Section 4.

**Theorem 2.4:** *Consider the wave equation with Kelvin-Voigt and viscous damping (1), (2), where $\sigma>0$, $\mu\geq 0$ are constants. Suppose that Assumption (H1) holds. Define for each $\omega>0$:*

$$r := \frac{\omega\sqrt{(\mu\sigma-1)^2\omega^2 + (\mu+\sigma\omega^2)^2}}{1+\sigma^2\omega^2}, \tag{20}$$

$$\begin{aligned} a &:= \sqrt{r}\cos(\theta/2) \\ b &:= \sqrt{r}\sin(\theta/2) \end{aligned} \tag{21}$$

*where $\theta\in(0,\pi)$ is the unique angle that satisfies the equations*

$$\cos(\theta) = \frac{(\mu\sigma-1)\omega^2}{(1+\sigma^2\omega^2)r}, \quad \sin(\theta) = \frac{(\mu+\sigma\omega^2)\omega}{(1+\sigma^2\omega^2)r}. \tag{22}$$

*Then the following inequality holds:*

$$\gamma_\infty \geq \sup_{\omega>0}\left(\sqrt{\frac{\max_{x\in[0,1]}\left(\cosh(2ax)-\cos(2bx)\right)}{\cosh(2a)-\cos(2b)}}\right). \tag{23}$$

Noticing that the right hand side of (23) is always greater or equal to 1, Theorem 2.3 allows us to obtain the following corollary for the sup norm.

**Corollary 2.5:** *Consider the wave equation with Kelvin-Voigt and viscous damping (1), (2), where $\sigma > 0$, $\mu \geq 0$ are constants. Suppose that $\mu\sigma = 1$. Then $\gamma_\infty = 1$.*

## 3. Asymptotic Gain in the L² Norm

*3.1. Existence of Asymptotic Gain in the L² Norm*

It is not clear whether Assumption (H1) holds for all $\sigma > 0$, $\mu \geq 0$ for the wave equation with Kelvin-Voigt and viscous damping. However, the analogue of Assumption (H1) for the $L^2$ norm holds for all $\sigma > 0$, $\mu \geq 0$. This is a consequence of the following theorem.

**Theorem 3.1:** *Consider the wave equation with Kelvin-Voigt and viscous damping (1), (2), where $\sigma > 0$, $\mu \geq 0$ are constants. There exists a constant $\gamma \geq 0$ such that for every bounded disturbance $d \in C^2(\Re_+)$ for which (1), (2) has a unique solution $u \in C^2(\Re_+ \times [0,1])$ with $(\partial u/\partial t)[t] \in C^2([0,1])$ for all $t \geq 0$, the following inequality holds:*

$$\limsup_{t \to +\infty}\left(\|u[t]\|_2\right) \leq \gamma \sup_{t \geq 0}\left(|d(t)|\right). \tag{24}$$

*Moreover, the constant $\gamma \geq 0$ satisfies the inequality*

$$\gamma \leq G(\mu, \sigma), \tag{25}$$

*where*

$$G(\mu, \sigma) := \begin{cases} \dfrac{1}{\sqrt{3}} & \text{if } \mu\sigma \geq 1 \\ \dfrac{1}{\pi}\sqrt{2\sum_{n=1}^{\infty} n^{-2} A_n^2} & \text{if } 0 \leq \mu\sigma < 1 \end{cases} \tag{26}$$

*and*

$$A_n = 1 + 2\sqrt{1-\mu\sigma}\left(\frac{2\beta_n\sqrt{1-\mu\sigma}}{\sigma(\mu+n^2\pi^2\sigma)(1+\beta_n) - 2\beta_n}\right)^{\beta_n}, \quad \beta_n := \frac{\mu + n^2\pi^2\sigma}{\sqrt{(\mu+n^2\pi^2\sigma)^2 - 4n^2\pi^2}},$$

$$\text{if } \mu + n^2\pi^2\sigma > 2n\pi \geq 2\sigma^{-1}, \tag{27}$$

$$A_n = 1, \text{ if } n\pi\sigma < 1 \text{ and } \mu + n^2\pi^2\sigma \geq 2n\pi, \tag{28}$$

$$A_n = 1 + 2\sqrt{1-\mu\sigma}\exp\left(-1 - \frac{1}{\sqrt{1-\mu\sigma}}\right), \text{ if } n\pi\sigma = 1 + \sqrt{1-\mu\sigma}, \tag{29}$$

$$A_n = 1 + 2\sqrt{1-\sigma\mu}\,\frac{\exp\left(\dfrac{(\mu+n^2\pi^2\sigma)}{\sqrt{4n^2\pi^2 - (\mu+n^2\pi^2\sigma)^2}}\arccos\left(\dfrac{2-\sigma\mu-n^2\pi^2\sigma^2}{2\sqrt{1-\sigma\mu}}\right)\right)}{\exp\left(\dfrac{(\mu+n^2\pi^2\sigma)\pi}{\sqrt{4n^2\pi^2 - (\mu+n^2\pi^2\sigma)^2}}\right) - 1}, \text{ if } \mu+n^2\pi^2\sigma < 2n\pi. \tag{30}$$

The proof of Theorem 3.1 can be found in the following section. A direct consequence of Theorem 3.1 is the fact that for every $\sigma > 0$, $\mu \geq 0$, we can define the asymptotic gain of the displacement in the $L^2$ norm for the wave equation with Kelvin-Voigt and viscous damping in the following way:

$$\gamma_2 := \sup \left\{ \frac{\limsup_{t \to +\infty}(\|u[t]\|_2)}{\sup_{t \geq 0}(|d(t)|)} : (f,g) \in X, d \in \Phi(f,g), d \neq 0 \right\}, \tag{31}$$

where $u \in C^2(\Re_+ \times [0,1])$ with $(\partial u/\partial t)[t] \in C^2([0,1])$ for all $t \geq 0$ is the solution of the initial-boundary value problem (1), (2), (16). Moreover, it follows from Theorem 3.1 that the following inequality holds:

$$\gamma_2 \leq G(\mu, \sigma), \tag{32}$$

where the function $G(\mu,\sigma)$ is defined by (26). Finally, if Assumption (H1) holds then the following inequality is a direct consequence of the fact that $\|u\|_2 \leq \|u\|_\infty$ for all $u \in C^0([0,1])$:

$$\gamma_2 \leq \gamma_\infty. \tag{33}$$

*3.2. Lower Bound of Asymptotic Gain in the $L^2$ Norm*

Our second main result in the present work concerning the $L^2$ norm is given below and it is proved by means of the frequency analysis procedure that was described above. Its proof is provided in Section 4.

**Theorem 3.2:** *Consider the wave equation with Kelvin-Voigt and viscous damping (1), (2), where $\sigma > 0$, $\mu \geq 0$ are constants. Define $r, a, b$ by means of (20), (21) for each $\omega > 0$, where $\theta \in (0, \pi)$ is the unique angle that satisfies equations (22). Then the following inequality holds:*

$$\gamma_2 \geq \sup_{\omega > 0}(Q(\omega)), \tag{34}$$

where $Q(\omega) := \sqrt{p + \sqrt{M}}$ with

$$M := \frac{4(a^2 + b^2) - 4(b\sin(2b)\cosh(2a) + a\sinh(2a)\cos(2b)) + \sinh^2(2a) + \sin^2(2b)}{16(a^2 + b^2)(\cosh(2a) - \cos(2b))^2}, \tag{35}$$

$$p = \frac{b\sinh(2a) - a\sin(2b)}{4ab(\cosh(2a) - \cos(2b))}. \tag{36}$$

Combining Theorem 3.1 and Theorem 3.2, we get the following corollary, which is proved in the following section.

**Corollary 3.3:** *Consider the wave equation with Kelvin-Voigt and viscous damping (1), (2), where $\sigma > 0$, $\mu \geq 0$ are constants. Suppose that $\mu\sigma \geq 1$. Then $\gamma_2 = 1/\sqrt{3}$.*

## 4. Proofs of Main Results

**Proof of Theorem 2.3:** The fact that Assumption (H1) holds when $2 < 2\mu\sigma + \sigma^2\pi^2$ is a direct consequence of Theorem 2.2 in [10]. We next notice that every solution $u \in C^2(\Re_+ \times [0,1])$ with $(\partial u/\partial t)[t] \in C^2([0,1])$ for all $t \geq 0$, of (1), (2) is a solution of the system of the integro-differential equation

$$\frac{\partial u}{\partial t}(t,x) = \sigma \frac{\partial^2 u}{\partial x^2}(t,x) - \frac{\mu\sigma - 1}{\sigma} u(t,x) + \frac{\mu\sigma - 1}{\sigma^2} \int_0^t \exp\left(-\frac{t-s}{\sigma}\right) u(s,x) ds \qquad (37)$$
$$+ \exp\left(-\frac{t}{\sigma}\right) \left( \frac{\partial u}{\partial t}(0,x) - \sigma \frac{\partial^2 u}{\partial x^2}(0,x) + \frac{\mu\sigma - 1}{\sigma} u(0,x) \right)$$

Let $\theta \in (0, \pi)$, $\varphi \in [0, \pi - \theta)$ be a pair of angles with $|\mu\sigma - 1| + 1 - \mu\sigma < \sigma^2\varphi^2$ (such a pair of angles exists due to the fact that $2 < 2\mu\sigma + \sigma^2\pi^2$). Define the positive function $\eta(x) := \sin(\theta + \varphi x)$ for $x \in [0,1]$ and the norm $\|u\|_{\infty,\eta} := \max_{0 \leq x \leq 1}(|u(x)|/\eta(x))$. Notice that Assumptions (H1), (H2), (H3), (H4) in [10] hold for the PDE problem

$$\frac{\partial u}{\partial t}(t,x) = \frac{\partial^2 u}{\partial x^2}(t,x) - \frac{\mu\sigma - 1}{\sigma} u(t,x) + f(t,x)$$

with $u(t,0) - d(t) = u(t,1) = 0$. Applying Corollary 4.3 in [10] to (37) with

$$f(t,x) = \frac{\mu\sigma - 1}{\sigma^2} \int_0^t \exp\left(-\frac{t-s}{\sigma}\right) u(s,x) ds + \exp\left(-\frac{t}{\sigma}\right) \left( \frac{\partial u}{\partial t}(0,x) - \sigma \frac{\partial^2 u}{\partial x^2}(0,x) + \frac{\mu\sigma - 1}{\sigma} u(0,x) \right)$$

and using the semigroup property and the fact that $\limsup_{t \to +\infty}(\|f[t]\|_{\infty,\eta}) \leq \frac{|\mu\sigma - 1|}{\sigma} \sup_{t \to +\infty}(\|u[t]\|_{\infty,\eta})$, we get

$$\limsup_{t \to +\infty}(\|u[t]\|_{\infty,\eta}) \leq \frac{1}{\sin(\theta)} \limsup_{t \to +\infty}(|d(t)|) + \frac{|\mu\sigma - 1|}{\mu\sigma - 1 + \sigma^2\varphi^2} \limsup_{t \to +\infty}(\|u[t]\|_{\infty,\eta})$$

It follows from the above inequality and the fact that $\|u\|_\infty \leq \|u\|_{\infty,\eta}$ for all $u \in C^0([0,1])$ that $\gamma_\infty \leq 1/\left( \left(1 - |s|/(s + \varphi^2)\right)\sin(\theta) \right)$ for all $\theta \in (0,\pi)$, $\varphi \in [0, \pi - \theta)$ with $|s| - s < \varphi^2$ and $s := (\mu\sigma - 1)/\sigma^2$. Therefore, for every $\theta \in (0,\pi)$ with $\theta < \pi - \sqrt{|s| - s}$, it holds that $\gamma_\infty \leq 1/\left( \left(1 - |s|/(s + (\pi - \theta)^2)\right)\sin(\theta) \right)$. Thus, inequality (19) holds with $g(s) := \inf\left\{ 1/(\sin(\theta)(1 - p(s,\theta))) : 0 < \theta < \pi - \sqrt{|s| - s} \right\}$ and $p(s,\theta) := |s|/(s + (\pi - \theta)^2)$. The proof is complete. ◁

**Proof of Theorem 2.4:** We apply the frequency analysis methodology for the parameterized family of inputs $d_\omega(t) = \sin(\omega t)$ with parameter $\omega > 0$. A periodic solution $u_\omega[t]$ of (1), (2) that corresponds to the input $d_\omega(t) = \sin(\omega t)$ is given by the formula

$$u_\omega(t,x) = \sin(\omega t)h(x) + \cos(\omega t)g(x), \text{ for } t \geq 0, x \in [0,1] \qquad (38)$$

where $h, g$ are solutions of the boundary-value problem

$$\sigma\omega g''(x) - \mu\omega g(x) - \omega^2 h(x) - h''(x) = 0$$
$$\omega^2 g(x) + g''(x) + \sigma\omega h''(x) - \mu\omega h(x) = 0,$$
(39)

$$h(0) = 1, \ h(1) = g(0) = g(1) = 0. \quad (40)$$

It can be verified that the functions

$$h(x) = \frac{\left(\cosh(a(2-x)) - \cos(2b)\cosh(ax)\right)\cos(bx) - \sin(2b)\sin(bx)\cosh(ax)}{\cosh(2a) - \cos(2b)},$$
$$g(x) = \frac{\sin(2b)\cos(bx)\sinh(ax) - \left(\sinh(a(2-x)) + \cos(2b)\sinh(ax)\right)\sin(bx)}{\cosh(2a) - \cos(2b)},$$
(41)

where $a, b$ are defined by (21), are solutions of the boundary-value problem (39), (40). Using (6) and (38), we obtain the inequality

$$\sup_{\omega>0}\left(\max_{0\le t\le 2\pi/\omega}\left(\max_{0\le x\le 1}\left(|\sin(\omega t)h(x) + \cos(\omega t)g(x)|\right)\right)\right) \le \gamma_\infty. \quad (42)$$

Since

$$\max_{0\le t\le 2\pi/\omega}\left(\max_{0\le x\le 1}\left(|\sin(\omega t)h(x) + \cos(\omega t)g(x)|\right)\right) = \max_{0\le x\le 1}\left(\max_{0\le t\le 2\pi/\omega}\left(|\sin(\omega t)h(x) + \cos(\omega t)g(x)|\right)\right) = \max_{0\le x\le 1}\left(\sqrt{h^2(x) + g^2(x)}\right) \quad (43)$$

and since (41) implies that

$$h^2(x) + g^2(x) = \frac{\cosh(2a(1-x)) - \cos(2b(1-x))}{\cosh(2a) - \cos(2b)}, \ \text{for } x \in [0,1], \quad (44)$$

we obtain from (42), (43) and (44) the desired inequality (23). The proof is complete. ◁

**Proof of Theorem 3.1:** It suffices to show that there exists a constant $\gamma \ge 0$ such that for every bounded disturbance $d \in C^2(\Re_+)$ for which (1), (2) has a unique solution $u \in C^2(\Re_+ \times [0,1])$ with $(\partial u/\partial t)[t] \in C^2([0,1])$ for all $t \ge 0$, with the following property:

"For every $\varepsilon > 0$ there exists $T > 0$ such that $\|u[t]\|_2 \le \gamma \sup_{s\ge 0}(|d(s)|) + \varepsilon$ for all $t \ge T$".

Let an arbitrary bounded disturbance $d \in C^2(\Re_+)$ for which (1), (2) has a unique solution $u \in C^2(\Re_+ \times [0,1])$ with $(\partial u/\partial t)[t] \in C^2([0,1])$ for all $t \ge 0$. Define:

$$y_n(t) = \sqrt{2}\int_0^1 u(t,x)\sin(n\pi x)dx, \ \text{for } n = 1, 2, \ldots. \quad (45)$$

Definition (45) in conjunction with (1), (2) implies that the following differential equations hold for all $t \ge 0$ and $n = 1, 2, \ldots$:

$$\ddot{y}_n(t) + \left(\mu + n^2\pi^2\sigma\right)\dot{y}_n(t) + n^2\pi^2 y_n(t) = n\pi\sqrt{2}\left(\sigma\dot{d}(t) + d(t)\right) \quad (46)$$

When $\mu + n^2\pi^2\sigma > 2n\pi$ then we may define

$$k_n := \frac{\mu + n^2\pi^2\sigma}{2}, \ r_n := \frac{\sqrt{\left(\mu + n^2\pi^2\sigma\right)^2 - 4n^2\pi^2}}{2} \quad (47)$$

and express the solution of (46) by the following formula:

$$y_n(t) = \frac{\exp(-(k_n - r_n)t)}{2r_n}\left((k_n + r_n - (k_n - r_n)\exp(-2r_n t))y_n(0) + \dot{y}_n(0)(1 - \exp(-2r_n t))\right) \quad (48)$$

$$-\frac{n\pi\sigma}{r_n\sqrt{2}}\exp(-(k_n - r_n)t)(1 - \exp(-2r_n t))d(0) + g_n(t),$$

$$g_n(t) := \frac{n\pi}{r_n\sqrt{2}}\int_0^t \left((\sigma k_n + \sigma r_n - 1)\exp(-(k_n + r_n)(t - \tau)) + (1 + \sigma r_n - \sigma k_n)\exp(-(k_n - r_n)(t - \tau))\right)d(\tau)d\tau. \quad (49)$$

When $n\pi\sigma = 1 \pm \sqrt{1 - \mu\sigma}$ (or equivalently when $\mu + n^2\pi^2\sigma = 2n\pi$) then the solution of (46) is given by the following formulas:

$$y_n(t) = (1 + k_n t)y_n(0)\exp(-k_n t) + \dot{y}_n(0)t\exp(-k_n t) - n\pi\sigma\sqrt{2}t\exp(-k_n t)d(0) + g_n(t), \quad (50)$$

$$g_n(t) = n\pi\sqrt{2}\int_0^t \left(\sigma + (1 - \sigma k_n)(t - s)\right)\exp(-k_n(t - s))d(s)ds. \quad (51)$$

When $\mu + n^2\pi^2\sigma < 2n\pi$ then we may define:

$$\omega_n := \frac{\sqrt{4n^2\pi^2 - \left(\mu + n^2\pi^2\sigma\right)^2}}{2} \quad (52)$$

and in this case the solution of (46) is given by the following formulas:

$$y_n(t) = y_n(0)\left(\frac{k_n}{\omega_n}\sin(\omega_n t) + \cos(\omega_n t)\right)\exp(-k_n t) + \frac{\dot{y}_n(0)}{\omega_n}\sin(\omega_n t)\exp(-k_n t) - \frac{n\pi\sigma\sqrt{2}}{\omega_n}\sin(\omega_n t)\exp(-k_n t)d(0) + g_n(t), \quad (53)$$

$$g_n(t) = \frac{n\pi\sqrt{2}}{\omega_n}\int_0^t \left(\sigma\omega_n \cos(\omega_n(t - s)) + (1 - \sigma k_n)\sin(\omega_n(t - s))\right)\exp(-k_n(t - s))d(s)ds. \quad (54)$$

We next show that for every $\theta > 0$ there exists $T > 0$ such that for all $t \geq T$ it holds that:

$$|y_n(t) - g_n(t)| \leq \frac{\theta}{n} \text{ for all } t \geq T, \ n = 1, 2, \ldots \quad (55)$$

Definitions (47) imply that there exist positive constants $\delta, \tilde{\beta}, c > 0$ so that the following inequalities hold for all $n \geq 1$ with $\mu + n^2\pi^2\sigma > 2n\pi$:

$$\delta \geq \frac{2n^2\pi^2}{\mu + n^2\pi^2\sigma} \geq k_n - r_n \geq \frac{n^2\pi^2}{\mu + n^2\pi^2\sigma} \geq \tilde{\beta}, \quad (56)$$

$$r_n \geq cn^2\pi^2.$$

It follows from (48) and (56) that the following estimate holds for all $t \geq 0$ and $n \geq 1$ with $\mu + n^2\pi^2\sigma > 2n\pi$:

$$|y_n(t) - g_n(t)| \leq \frac{\exp(-\tilde{\beta}t)}{2r_n}\left((k_n + r_n)|y_n(0)| + |\dot{y}_n(0)|\right) + \frac{\sigma}{cn\pi\sqrt{2}}\exp(-\tilde{\beta}t)|d(0)|. \quad (57)$$

Definition (45) and (2) imply that $y_n(t) = \sqrt{2}d(t)/(n\pi) + (\sqrt{2}/(n\pi))\int_0^1 u'(t, x)\cos(n\pi x)dx$. Therefore, we get:

$$|y_n(t)| \leq \sqrt{2}\frac{|d(t)|}{n\pi} + \frac{1}{n\pi}\|u'[t]\|_2, \text{ for all } t \geq 0, \ n = 1, 2, \ldots \quad (58)$$

Similarly, definition (32) implies that:

$$|\dot{y}_n(t)| \leq \left\|\frac{\partial u}{\partial t}[t]\right\|_2, \text{ for all } t \geq 0 \text{ and } n \geq 1. \tag{59}$$

Using (56), (57), (58) and (59), we obtain for all $t \geq 0$ and $n \geq 1$ with $\mu + n^2\pi^2\sigma > 2n\pi$:

$$|y_n(t) - g_n(t)| \leq \frac{\exp(-\tilde{\beta}t)}{2r_n}\left((\delta + 2r_n)\frac{1}{n\pi}\left(\sqrt{2}|d(0)| + \|u'[0]\|_2\right) + \left\|\frac{\partial u}{\partial t}[0]\right\|_2\right) + \frac{\sigma}{cn\pi\sqrt{2}}\exp(-\tilde{\beta}t)|d(0)|. \tag{60}$$

It follows from (56) and (60) that the following inequality holds for all $t \geq 0$ and $n \geq 1$ with $\mu + n^2\pi^2\sigma > 2n\pi$:

$$|y_n(t) - g_n(t)| \leq \frac{\exp(-\tilde{\beta}t)}{2n\pi}\left(\left(2 + \frac{\delta}{c\pi^2}\right)\left(\sqrt{2}|d(0)| + \|u'[0]\|_2\right) + \frac{1}{c\pi}\left\|\frac{\partial u}{\partial t}[0]\right\|_2 + \frac{\sigma}{c}\sqrt{2}|d(0)|\right). \tag{61}$$

If $\mu + n^2\pi^2\sigma > 2n\pi$ holds for all integers $n \geq 1$ then (55) is a direct consequence of (61). We next consider the case that there exist integers $n \geq 1$ for which $\mu + n^2\pi^2\sigma \leq 2n\pi$.

Let $N \geq 1$ be the integer defined in the following way: $N = 1$ if $\mu\sigma > 1$ and $N = \left[\left(1 + \sqrt{1 - \mu\sigma}\right)/(\sigma\pi)\right] + 1$ if $\mu\sigma \leq 1$. Notice that $\mu + n^2\pi^2\sigma > 2n\pi$ for all $n \geq N$. Without loss of generality, we may assume that the positive constant $\tilde{\beta} > 0$ for which (56) holds, also satisfies the inequality $(\mu + \pi^2\sigma)/2 > \tilde{\beta}$. It then follows from (47), (50) and (53) that there exists a constant $\tilde{G} > 0$ (depending on the initial conditions and $d(0)$) such that the following inequality holds for all $t \geq 0$ and $n = 1, ..., N - 1$ with $\mu + n^2\pi^2\sigma \leq 2n\pi$:

$$|y_n(t) - g_n(t)| \leq \tilde{G}\exp(-\tilde{\beta}t).$$

The above inequality implies that

$$|y_n(t) - g_n(t)| \leq G\frac{N-1}{n}\exp(-\beta t). \tag{62}$$

Therefore, if there exist integers $n \geq 1$ for which $\mu + n^2\pi^2\sigma \leq 2n\pi$, then (55) is a direct consequence of (61) and (62). Thus, we have shown that (55) holds.

For all $n \geq 1$ with $\mu + n^2\pi^2\sigma > 2n\pi$, definitions (47), (49) imply the following estimates for all $t \geq 0$:

$$|g_n(t)| \leq \frac{n\pi}{r_n\sqrt{2}}\|d\|_\infty \int_0^t |(\sigma k_n + \sigma r_n - 1)\exp(-(k_n + r_n)(t - \tau)) + (1 + \sigma r_n - \sigma k_n)\exp(-(k_n - r_n)(t - \tau))|d\tau$$

$$\leq \frac{n\pi}{r_n\sqrt{2}}\|d\|_\infty\left(\frac{1}{k_n - r_n} - \frac{1}{k_n + r_n}\right)$$

if $\mu\sigma \geq 1$ or $n\pi\sigma \leq 1$ (equivalently if $\sigma(k_n - r_n) \leq 1$) \hfill (63)

$$|g_n(t)| \leq \frac{n\pi}{r_n\sqrt{2}}\|d\|_\infty \int_0^t |(\sigma k_n + \sigma r_n - 1)\exp(-(k_n + r_n)(t - \tau)) + (1 + \sigma r_n - \sigma k_n)\exp(-(k_n - r_n)(t - \tau))|d\tau$$

$$\leq \frac{n\pi}{r_n\sqrt{2}}\|d\|_\infty\left(\frac{1}{k_n - r_n} - \frac{1}{k_n + r_n}\right)\left(1 + 2\left(\frac{\sigma k_n - \sigma r_n - 1}{\sigma k_n + \sigma r_n - 1}\right)^{\frac{k_n}{2r_n}}\sqrt{(\sigma k_n + \sigma r_n - 1)(\sigma k_n - \sigma r_n - 1)}\right)$$

if $\mu\sigma < 1$ and $n\pi\sigma > 1$ (equivalently if $\sigma(k_n - r_n) > 1$), \hfill (64)

where $\|d\|_\infty := \sup_{s \geq 0}(|d(s)|)$. Notice that $\mu\sigma \geq 1$ (which implies $\mu + n^2\pi^2\sigma \geq 2n\pi$ for all $n \geq 1$), definitions (47), (51) and estimate (63) give:

$$|g_n(t)| \leq \frac{\sqrt{2}}{n\pi}\|d\|_\infty, \text{ for all } n \geq 1 \text{ and } t \geq 0. \tag{65}$$

Using (63), (64), (51), (54), (47), (52) and definitions (27), (28), (30), we obtain the following estimates for the case $\mu\sigma < 1$:

$$|g_n(t)| \leq \frac{\sqrt{2}}{n\pi} A_n \|d\|_\infty, \text{ for all } n \geq 1 \text{ and } t \geq 0. \tag{66}$$

We are ready to finish the proof. Let arbitrary $\varepsilon > 0$ be given. Since $\lim_{n \to +\infty}(A_n) = 1$ (a consequence of definition (27)), it follows that the sequence $\{A_n\}_{n=1}^\infty$ is bounded. Let $\bar{A}$ be an upper bound for the sequence $\{A_n\}_{n=1}^\infty$. Pick $\theta > 0$ sufficiently small so that $\theta^2(\pi^2/6) + \theta(\pi\sqrt{2}/3)\|d\|_\infty \bar{A} \leq \varepsilon^2$. Then there exists $T > 0$ so that (55) holds. Using the triangle inequality, we obtain from (55) for all $n \geq 1$ and $t \geq T$:

$$|y_n(t)|^2 \leq \frac{\theta^2}{n^2} + |g_n(t)|^2 + \frac{2\theta}{n}|g_n(t)|. \tag{67}$$

Definition (45) and the fact that the set of functions $\{\sqrt{2}\sin(n\pi x): n = 1,2,...\}$ is an orthonormal basis of $L^2(0,1)$ implies that Parseval's identity holds, i.e., $\|u[t]\|_2^2 = \sum_{n=1}^\infty y_n^2(t)$. If $\mu\sigma < 1$ then we get from (66), (67) for all $t \geq T$:

$$\|u[t]\|_2^2 \leq \theta^2 \sum_{n=1}^\infty n^{-2} + \frac{2}{\pi^2}\|d\|_\infty^2 \sum_{n=1}^\infty n^{-2} A_n^2 + \frac{2\sqrt{2}}{\pi}\theta\|d\|_\infty \sum_{n=1}^\infty n^{-2} A_n. \tag{68}$$

Using (68) and the fact $\sum_{n=1}^\infty n^{-2} = \pi^2/6$, we get for all $t \geq T$:

$$\|u[t]\|_2^2 \leq \theta^2 \frac{\pi^2}{6} + (G(\mu,\sigma))^2 \|d\|_\infty^2 + \frac{\pi\sqrt{2}}{3}\theta\|d\|_\infty \bar{A}, \tag{69}$$

where $\bar{A}$ is the upper bound for the sequence $\{A_n\}_{n=1}^\infty$. It follows from (69) and the fact that $\theta^2(\pi^2/6) + \theta(\pi\sqrt{2}/3)\|d\|_\infty \bar{A} \leq \varepsilon^2$ that $\|u[t]\|_2 \leq G(\mu,\sigma)\|d\|_\infty + \varepsilon$ for all $t \geq T$.

The analysis is similar for the case $\mu\sigma \geq 1$ (simply set $A_n \equiv 1$ in the above analysis).
The proof is complete. ◁

**Proof of Theorem 3.2:** We apply the frequency analysis methodology for the parameterized family of inputs $d_\omega(t) = \sin(\omega t)$ with parameter $\omega > 0$. A periodic solution $u_\omega[t]$ of (1), (2) that corresponds to the input $d_\omega(t) = \sin(\omega t)$ is given by (38), where $h, g$ are given by (41). Using (38), we get for all $t \geq 0$:

$$\|u[t]\|_2^2 = p + q_1 \cos(2\omega t) + q_2 \sin(2\omega t), \tag{70}$$

where

$$p := \frac{1}{2}\int_0^1 (h^2(x) + g^2(x))dx, \tag{71}$$

$$q_1 := \frac{1}{2}\int_0^1 (g^2(x) - h^2(x))dx, \tag{72}$$

$$q_2 := \int_0^1 h(x)g(x)dx. \tag{73}$$

Using (41), (44), (71), (72), (73), we obtain (36) as well as the following formulas:

$$q_1 = \frac{\cosh(2a)\cos(2b)-1}{2(\cosh(2a)-\cos(2b))^2} + \frac{b\sin(2b)-a\sinh(2a)}{4(a^2+b^2)(\cosh(2a)-\cos(2b))} \tag{74}$$

$$q_2 = \frac{\sinh(2a)\sin(2b)}{2(\cosh(2a)-\cos(2b))^2} - \frac{a\sin(2b)+b\sinh(2a)}{4(a^2+b^2)(\cosh(2a)-\cos(2b))} \tag{75}$$

Equation (70) implies that

$$\max\left\{\|u[t]\|_2^2 : 0 \le t \le \frac{2\pi}{\omega}\right\} = p + \sqrt{q_1^2 + q_2^2}, \tag{76}$$

Equation (4.41) combined with (74), (75) and (6) gives (34) with $Q(\omega) := \sqrt{p+\sqrt{M}}$ and $M$ given by (35). The proof is complete. ◁

**Proof of Corollary 3.3:** It suffices to show that $\lim_{\omega \to 0^+}(Q(\omega)) = 1/\sqrt{3}$, where $Q(\omega) := \sqrt{p+\sqrt{M}}$ with $p, M$ defined by (35), (36). Indeed, definitions (35), (36) imply that $p \to 1/6$ and $M \to 1/36$ as $(a,b) \to (0,0)$. Moreover, definitions (20), (21) imply that $a \to 0$ and $b \to 0$ as $\omega \to 0^+$. The proof is complete. ◁

## 5. Graphical Illustration of the Theorem Statements

This section is devoted to the graphical illustration of the theorems' statements for the upper and lower bounds of the asymptotic gains. To this purpose, we define:

$$U_\infty := g\left(\frac{\mu\sigma-1}{\sigma^2}\right), \text{ for } \sigma > 0, \mu \ge 0 \text{ with } 2 < 2\mu\sigma + \sigma^2\pi^2, \tag{77}$$

where $g(s) := \inf\left\{1/(\sin(\theta)(1-p(s,\theta))) : 0 < \theta < \pi - \sqrt{|s|-s}\right\}$ and $p(s,\theta) := |s|/(s+(\pi-\theta)^2)$,

$$L_\infty := \sup_{\omega > 0}(\bar{A}(\omega)), \tag{78}$$

$$\bar{A}(\omega) := \sqrt{\frac{\max_{x\in[0,1]}(\cosh(2ax)-\cos(2bx))}{\cosh(2a)-\cos(2b)}}, \tag{79}$$

and $a, b$ are being defined by (21),

$$L_2 := \sup_{\omega > 0}(Q(\omega)), \tag{80}$$

where

$$Q(\omega) := \sqrt{p+\sqrt{M}}, \tag{81}$$

and $p, M$ are being defined by (35), (36) and

$$U_2 := G(\mu, \sigma), \tag{82}$$

where $G(\mu, \sigma)$ is defined by (26). Notice that the results of the previous sections guarantee that

$$\gamma_\infty \le U_\infty, \text{ for all } \sigma > 0, \mu \ge 0 \text{ with } 2 < 2\mu\sigma + \sigma^2\pi^2, \tag{83}$$

$$L_\infty \le \gamma_\infty, \text{ for all } \sigma > 0, \mu \ge 0 \text{ for which Assumption (H1) holds,} \tag{84}$$

$$L_2 \leq \gamma_2 \leq U_2, \text{ for all } \sigma > 0, \mu \geq 0. \tag{85}$$

Fig. 1 and Fig. 3 depict the lower and upper bounds of the asymptotic gain in the sup norm for $\sigma = 1$ and $\sigma = 0.5$, respectively, and for a wide range of values for $\mu$.

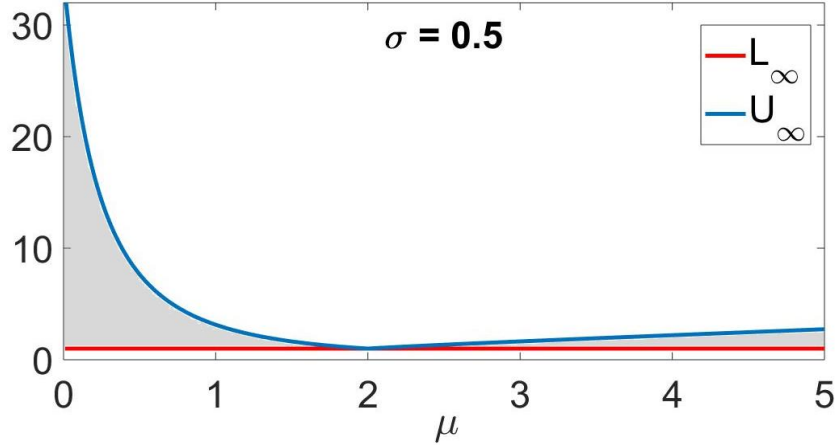

**Fig. 3:** The lower (red) and upper (blue) bounds for the asymptotic gain in the sup norm for $\sigma = 0.5$. The grey area depicts possible values for the asymptotic gain in the sup norm.

Fig. 4, Fig. 5 and Fig. 6 are Bode-like plots for the logarithm of $\bar{A}(\omega)$ defined by (79) for $\sigma = 0.0001$, $\sigma = 0.001$ and $\sigma = 0.01$, respectively, and for four different values of $\mu$. These plots indicate that for small $\sigma$ and $\mu$, $\bar{A}(\omega)$ presents "spikes" at frequencies which differ by $\pi$. On the other hand, for sufficiently large $\sigma$ and $\mu$, $\bar{A}(\omega)$ identically equal to 1.

Fig. 2 and Fig. 7 depict the lower and upper bounds of the asymptotic gain in the $L^2$ norm, for $\sigma = 1$ and $\sigma = 0.5$, respectively, and for a wide range of $\mu$ values. Notice that for $\mu\sigma \geq 1$, we get $L_2 = U_2 = 1/\sqrt{3}$. Comparing Fig. 2 and Fig. 7 with Fig. 1 and Fig. 3, we conclude that the asymptotic gain in the $L^2$ norm is estimated much more accurately than the asymptotic gain in the sup norm.

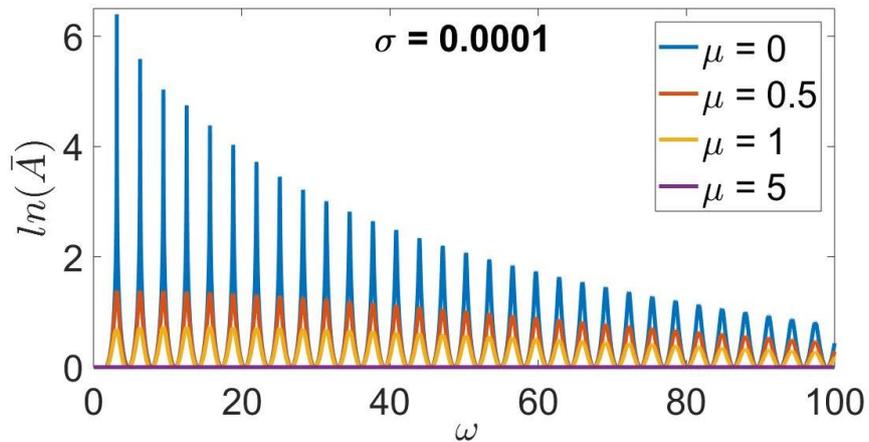

**Fig. 4:** Bode-like plot of $\ln(\bar{A}(\omega))$ for $\sigma = 0.0001$.

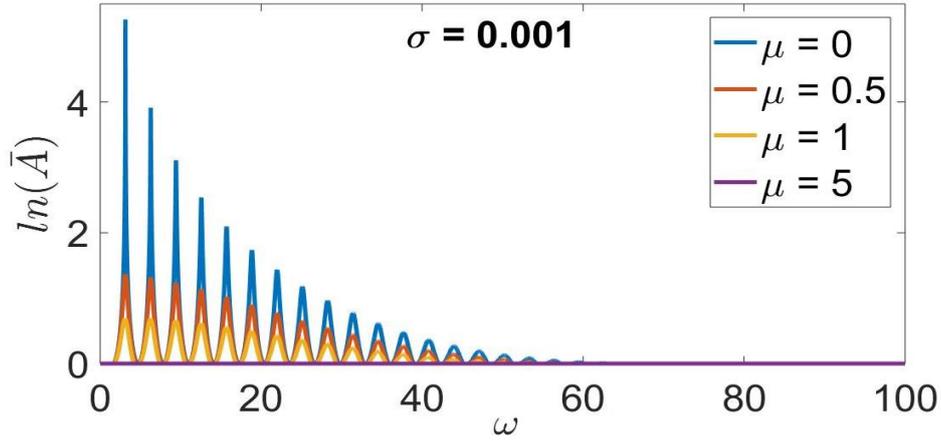

**Fig. 5:** Bode-like plot of $\ln(\bar{A}(\omega))$ for $\sigma = 0.001$.

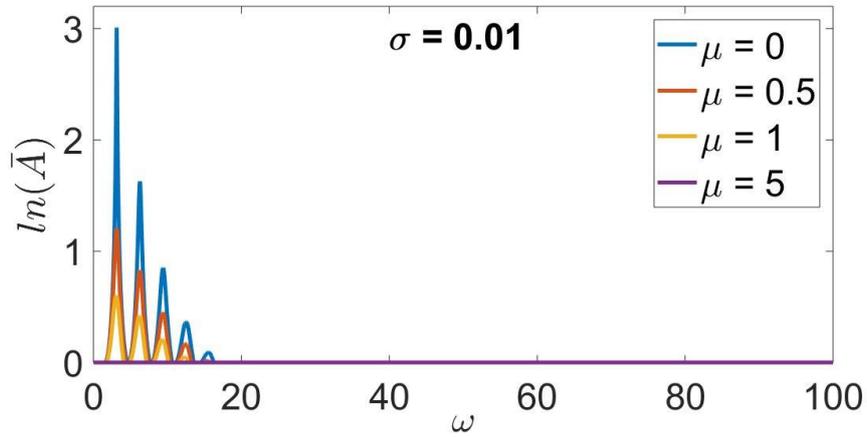

**Fig. 6:** Bode-like plot of $\ln(\bar{A}(\omega))$ for $\sigma = 0.01$.

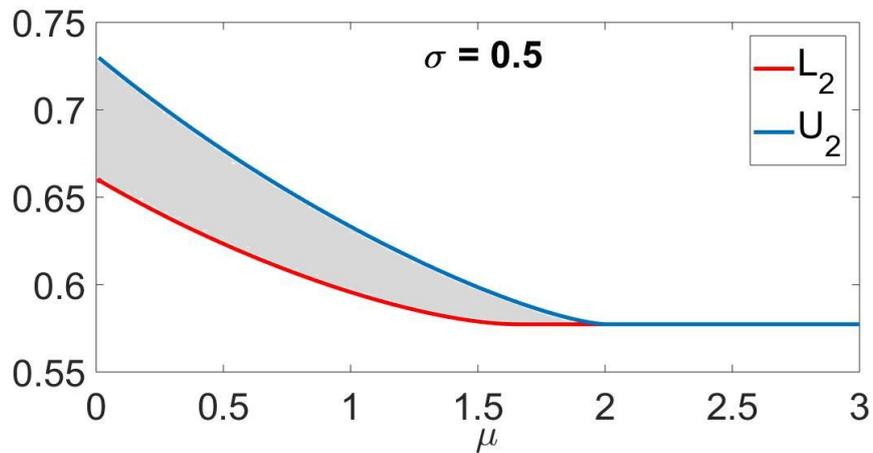

**Fig. 7:** The lower (red) and upper (blue) bounds for the asymptotic gain in the $L^2$ norm for $\sigma = 0.5$. The grey area depicts possible values for the asymptotic gain in the $L^2$ norm.

Fig. 8, Fig. 9, Fig. 10 and Fig. 11 are Bode-like plots for the logarithm of $Q(\omega)$ defined by (81) for $\sigma = 0.0001$, $\sigma = 0.001$, $\sigma = 0.01$ and $\sigma = 0.5$, respectively, and for four different values of $\mu$. As in the sup norm, these plots indicate that for small $\sigma$ and $\mu$, $Q(\omega)$ presents "spikes" at frequencies which differ by $\pi$. These are "dangerous" frequencies: the amplitude of the oscillation of the boundary of a string becomes magnified in its domain.

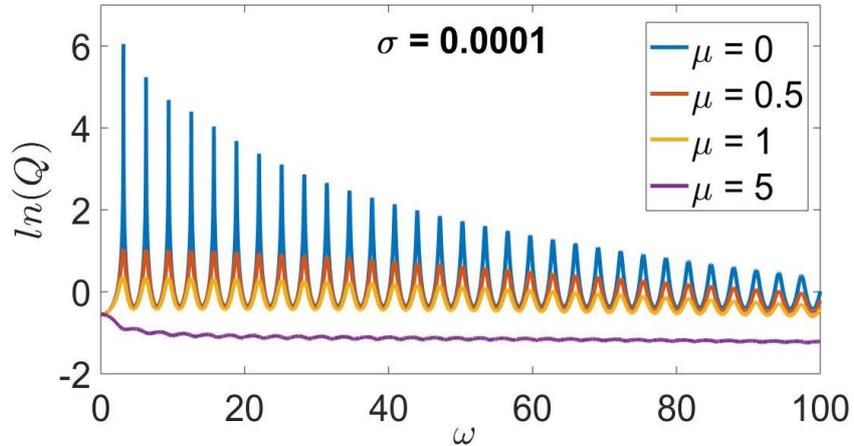

**Fig. 8:** Bode-like plot of $\ln(Q(\omega))$ for $\sigma = 0.0001$.

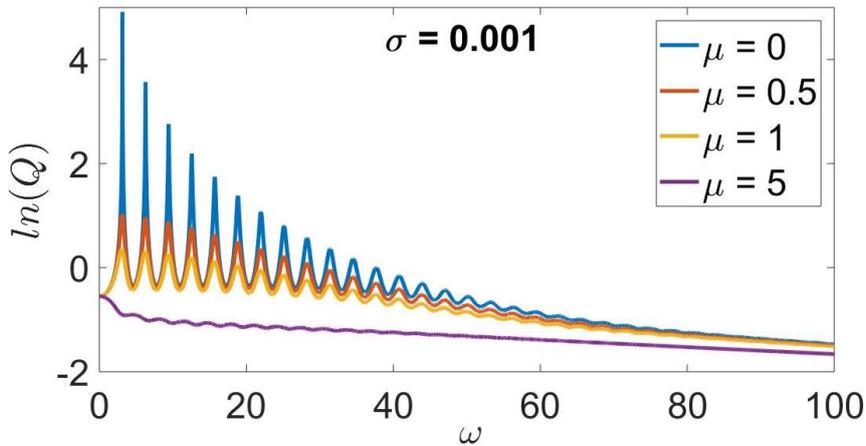

**Fig. 9:** Bode-like plot of $\ln(Q(\omega))$ for $\sigma = 0.001$.

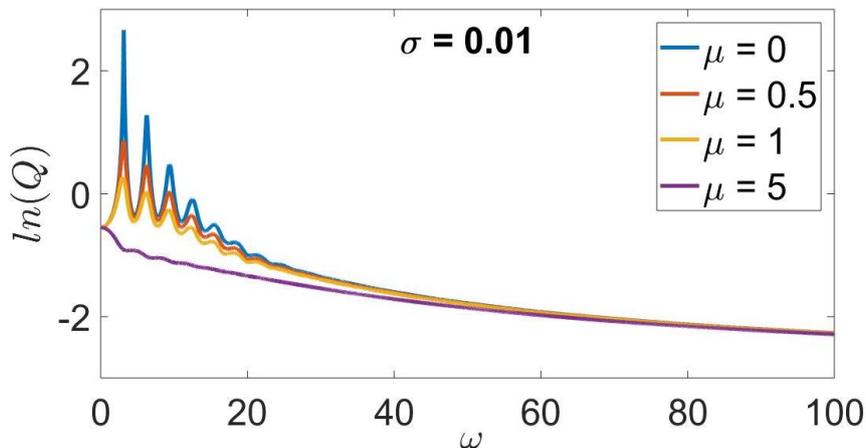

**Fig. 10:** Bode-like plot of $\ln(Q(\omega))$ for $\sigma = 0.01$.

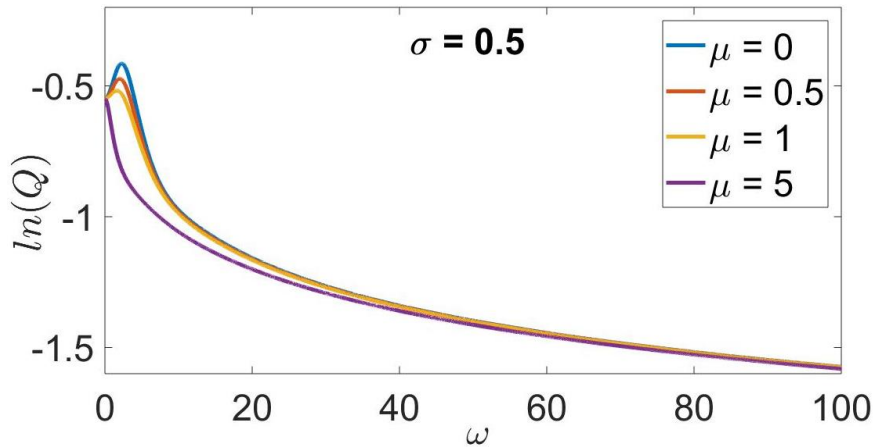

**Fig. 11:** Bode-like plot of $\ln(Q(\omega))$ for $\sigma = 0.5$.

## 6. Concluding Remarks

We have provided estimates for the asymptotic gains of the displacement of a vibrating string with external forcing. By considering an external boundary disturbance for the wave equation with Kelvin-Voigt and viscous damping, we have shown that the asymptotic gain property holds in the spatial $L^2$ norm of the displacement without any assumption for the damping coefficients. We have also provided upper and lower bounds for the asymptotic gains in the spatial sup-norm and the spatial $L^2$ norm of the displacement.

As noted in the introduction, this is the first systematic study of the asymptotic gains for a PDE with a boundary disturbance. Although here we studied the wave equation with Kelvin-Voigt and viscous damping, it is expected that similar studies will be performed for many other important PDEs.